\documentclass[12pt,reqno]{amsart}%
\usepackage{amsfonts}
\usepackage{amsmath}
\usepackage{amssymb}
\usepackage{graphicx}%
\providecommand{\U}[1]{\protect\rule{.1in}{.1in}}
\theoremstyle{plain}

\numberwithin{equation}{section}
\begin{document}
\title[Best constants and maximizers   for weighted Moser-Trudinger inequalities ]{Best constants and existence of maximizers   for weighted Moser-Trudinger inequalities}
\author{Mengxia Dong}
\address{Department of Mathematics\\
Wayne State University\\
Detroit, MI 48202, USA}
\email{mengxia.dong@wayne.edu}
\author{Guozhen Lu}
\address{Department of Mathematics\\
Wayne State University\\
Detroit, MI 48202, USA}
\email{gzlu@wayne.edu}

\thanks{Research of this work was partly supported by a US NSF grant DMS-1301595.}

\begin{abstract}

Sharp Moser-Trudinger inequalities and existence of maximizers for such inequalities   play an important role in geometric analysis, partial differential equations and other branches of modern mathematics. Such geometric inequalities have been studied extensively by many authors in recent years and there is a vast literature.
In this paper, we will establish the best constants for certain classes of weighted Moser-Trudinger inequalities on the entire Euclidean spaces $\mathbb{R}^N$. More precisely,
for given $N\ge 2$, \, $-\infty<s\le t<N$, \, $0<\alpha<\alpha_{N,t}:=(N-t)\omega^{\frac{1}{N-1}}_{N-1}$ and $\Phi_N(t):=\sum^{\infty}_{j=N-1}\frac{t^j}{j!}$,  we will show
  there exists a positive constant $C=C(N,s,t,\alpha)$ such that the following two inequalities
 \begin{equation}\label{equ0.1}
\int_{\mathbb{R}^N}\Phi_N(\alpha|u|^{N/(N-1)})\frac{dx}{|x|^t}\le C(\int_{\mathbb{R}^N}\frac{|u|^{N}}{|x|^s}dx)^{\frac{N-t}{N-s}}
\end{equation}
\begin{equation}\label{equ0.2}
\int_{\mathbb{R}^N}e^{\alpha|u|^{N/(N-1)}}|u|^N\frac{dx}{|x|^t}\le C(\int_{\mathbb{R}^N}\frac{|u|^{N}}{|x|^s}dx)^{\frac{N-t}{N-s}}
\end{equation}
holds for all   functions $u\in L^N(\mathbb{R}^N;|x|^{-s}dx)\cap \dot{W}^{1,N}(\mathbb{R}^N)$ with $\|\nabla u\|_{L^N(\mathbb{R}^N)}\le 1$.
 Moreover, the constant $\alpha_{N,t}=(N-t)\omega^{\frac{1}{N-1}}_{N-1}$ is sharp in the sense if $\alpha>\alpha_{N,t}$, then none of the above inequalities can hold with a uniform constant $C$ for all such $u$.  We will also prove the existence of maximizers of these sharp weighted inequalities. The class of functions considered here are not necessarily spherically symmetric. Our inequality (\ref{equ0.1}) (Theorem 1.1)  improves the earlier one
where such type of inequality was only considered for spherically symmetric functions by M. Ishiwata, M. Nakamura, H. Wadade in \cite{INW} (except in the case $s\not=0$).   Since
$\int_{\mathbb{R}^N}\Phi_N(\alpha|u|^{N/(N-1)})\frac{dx}{|x|^t}\le \int_{\mathbb{R}^N}e^{\alpha|u|^{N/(N-1)}}|u|^N\frac{dx}{|x|^t}$, our inequality
(\ref{equ0.2}) is stronger than inequality (\ref{equ0.1}). We can also replace the weight $\frac{1}{|x|^N}$ in Theorem 1.2 by  $\frac{1}{|x|^q}$ for $q>N$ (Theorem 1.3). 

We note that it suffices for us to prove the above inequalities for all functions not necessarily radially symmetric when $s=t$ by the well-known Caffareli-Kokn-Nirenberg inequalities \cite{CKN}.

\end{abstract}
\maketitle

\section{Introduction and main results}

In this article, our main purpose is to establish the weighted Moser-Trudinger type inequalities with   sharp constants and consider the existence of a maximizer associated with the weighted Moser-Trudinger type inequalities. The method we develop here does not need to assume the functions under consideration to be radially symmetric.\\

It is well known that Sobolev embedding gives us continuous embedding $W^{k,p}_0(\Omega)\subset L^q(\Omega)$ where $kp<N$ when $\Omega\subset \mathbb{R}^N (N\ge2)$ is a bounded domain with $1\le q\le\frac{Np}{N-kp}$, and $W^{k,p}_0(\Omega)\subset L^q(\Omega)$ for $1\le q<\infty$. However, it is not hard to show that in general $W^{1,N}_0(\Omega)\nsubseteq L^{\infty}(\Omega)$. In this case, Yudovich \cite{Yu}, Pohozaev \cite{Po} and Trudinger \cite{Tr} proved independently that $W^{1, N}_0(\Omega)\subset L_{\varphi_N}(\Omega)$,  where $L_{\varphi_N}(\Omega)$ is the Orlicz space associated with the Young function $\varphi_N(t)=\exp(\beta|t|^{N/N-1})-1$ for some $\beta>0$.  J. Moser proved the following sharp result in his 1971 paper \cite{Mo}:\\
\\
\textbf{Theorem A.} \emph{Let $\Omega$ be a domain with finite measure in Euclidean N-space $\mathbb{R}^N$, $N\ge2$. Then there exists a sharp constant $\alpha_{N}=N\left(
\frac{N\pi^{\frac{N}{2}}}{\Gamma(\frac{N}{2}+1)}\right)  ^{\frac{1}{N-1}}%
$ such that
\begin{displaymath}
\frac{1}{|\Omega|}\int_{\Omega}\exp(\beta|u|^{\frac{N}{N-1}})dx\le c_0
\end{displaymath}
for any $\beta\le\alpha_N$, any $u\in W^{1,N}_0(\Omega)$ with $\int_{\Omega}|\nabla u|^N dx\le1$. This constant $\beta_N$ is sharp in the sense that if $\beta>\beta_N$, then the above inequality can no longer hold with some $c_0$ independent of $u$}.\\

There are many generalizations related to the above classical Moser-Trudinger inequality, in particular to unbounded domains. As a scaling invariant form in $\mathbb{R}^N$, Adachi and Tanaka \cite{AT} proved the following inequality on the entire Euclidean space $\mathbb{R}^N$:\\
\\
\textbf{Theorem B.} \emph{For $N\ge2$ and $0<\alpha<\alpha_N=N\left(
\frac{N\pi^{\frac{N}{2}}}{\Gamma(\frac{N}{2}+1)}\right)  ^{\frac{1}{N-1}}%
$,  there exists a positive constant $C=C(N,\alpha)$ such that the inequality
\begin{displaymath}
\int_{\mathbb{R}^N}\Phi_N(\alpha|u(x)|^{N'})dx\le C\|u\|^N_{L^N(\mathbb{R}^N)}
\end{displaymath}
holds for all $u\in W^{1,N}(\mathbb{R}^N)$ with $\|\nabla u\|_{L^N(\mathbb{R}^N)}\le1$, where
\begin{displaymath}
\Phi_N(t):=\sum^{\infty}_{j=N-1}\frac{t^j}{j!}, \quad\quad t\ge0.
\end{displaymath}
Moreover, the constant $\alpha_N$ is sharp in the sense that    if $\alpha\ge\alpha_N$ then the inequality cannot hold with a uniform constant $C$ independent of $u$.}\\

Recently, Lam, Lu and Zhang proved in \cite{LamLuZhang} the precise asymptotic estimates for the following supremum.\\
\\
 \textbf{Theorem C.} \textit{Let }$N\geq2$\textit{, }$\alpha_{N}=N\left(
\frac{N\pi^{\frac{N}{2}}}{\Gamma(\frac{N}{2}+1)}\right)  ^{\frac{1}{N-1}}%
$,$~0\leq\beta<N$ \textit{and }$0\leq\alpha<\alpha_{N}.$ Denote%
\[
AT\left(  \alpha,\beta\right)  =\sup_{\left\Vert \nabla u\right\Vert _{N}%
\leq1}\frac{1}{\left\Vert u\right\Vert _{N}^{N-\beta}}\int_{%
\mathbb{R}
^{N}}\phi_{N}\left(  \alpha\left(  1-\frac{\beta}{N}\right)  \left\vert
u\right\vert ^{\frac{N}{N-1}}\right)  \frac{dx}{\left\vert x\right\vert
^{\beta}}.
\]
Then there exist positive constants $c=c\left(  N,\beta\right)  $ and
$C=C\left(  N,\beta\right)  $ such that when $\alpha$ is close enough to
$\alpha_{N}:$
\begin{equation}
\frac{c\left(  N,\beta\right)  }{\left(  1-\left(  \frac{\alpha}{\alpha_{N}%
}\right)  ^{N-1}\right)  ^{\left(  N-\beta\right)  /N}}\leq AT\left(
\alpha,\beta\right)  \leq\frac{C\left(  N,\beta\right)  }{\left(  1-\left(
\frac{\alpha}{\alpha_{N}}\right)  ^{N-1}\right)  ^{\left(  N-\beta\right)
/N}}. \label{1.3.1}%
\end{equation}
Moreover, the constant $\alpha_{N}$ is sharp in the sense that $AT\left(
\alpha_{N},\beta\right)  =\infty.$

The upper bound in the above estimates for the subcritical case was obtained by an argument inspired by the work of Lam and the second author \cite{LamLu} where 
a local Trudinger-Moser inequality on the level sets of the functions under consideration can lead to a global one on the entire spaces, without a priori knowing the validity of the critical inequality.

We remark that in dimension two, the upper bound for the $AT(\alpha, \beta)$ was also obtained in \cite{CST} using the critical Trudinger-Moser inequality in \cite{R}.  
We also note in the above theorem, we only impose the restriction on the norm
$\int_{\mathbb{R}^{N}}\left\vert \nabla u\right\vert ^{N}$ without restricting
the full norm
\[
\left[  \int_{\mathbb{R}^{N}}\left\vert \nabla u\right\vert ^{N}+\tau
\int_{\mathbb{R}^{N}}\left\vert u\right\vert ^{N}\right]  ^{1/N}\leq1.
\]
The method in \cite{AT} requires a symmetrization argument which is not
available in many other non-Euclidean settings. The above inequality fails at
the critical case $\alpha=\alpha_{N}$. So it is natural to ask when the above
can be true when $\alpha=\alpha_{N}$. This is done in Ruf \cite{R} and Li-Ruf \cite{LR} by
using the restriction on the full norm  the Sobolev space
$W^{1,N}\left(  \mathbb{R}^{N}\right)$.  \\
\\
\textbf{Theorem D.} \textit{For all
}$0\leq\alpha\leq\alpha_{N}:$%
\begin{equation}
\underset{\left\Vert u\right\Vert \leq1}{\sup}\int_{%
\mathbb{R}
^{N}}\phi_{N}\left(  \alpha\left\vert u\right\vert ^{\frac{N}{N-1}}\right)
dx<\infty\label{1.4}%
\end{equation}
\textit{where }%
\[
\left\Vert u\right\Vert =\left(  \int_{%
\mathbb{R}
^{N}}\left(  \left\vert \nabla u\right\vert ^{N}+\left\vert u\right\vert
^{N}\right)  dx\right)  ^{1/N}.
\]
\textit{Moreover, this constant }$\alpha_{N}$\textit{ is sharp in the sense
that if }$\alpha>\alpha_{N}$\textit{, then the supremum is infinity.}  \\

Surprisingly, Lam, Lu and Zhang have shown in \cite{LamLuZhang} that the subcritical Moser-Trudinger inequality in \cite{AT} and the critical Moser-Trudinger inequality in \cite{R, LR} are actually equivalent. Then we will provide another proof of the sharp critical Moser-Trudinger
inequality using the subcritical one, and vice versa. Furthermore, we have shown the following precise relationship between the supremums in the critical and subcritical Moser-Trudinger inequalities.\\
\\
\textbf{Theorem E.} \textit{Let }$N\geq2$\textit{,}$~0\leq\beta<N,~0<a,~b.$ Denote
\begin{align*}
MT_{a,b}\left(  \beta\right)   &  =\sup_{\left\Vert \nabla u\right\Vert
_{N}^{a}+\left\Vert u\right\Vert _{N}^{b}\leq1}\int_{%
\mathbb{R}
^{N}}\phi_{N}\left(  \alpha_{N}\left(  1-\frac{\beta}{N}\right)  \left\vert
u\right\vert ^{\frac{N}{N-1}}\right)  \frac{dx}{\left\vert x\right\vert
^{\beta}};\\
MT\left(  \beta\right)   &  =MT_{N,N}\left(  \beta\right)  .
\end{align*}
Then $MT_{a,b}\left(  \beta\right)  <\infty$ if and only if $b\leq N$. The
constant $\alpha_{N}$ is sharp. Moreover, we have the following identity:%
\begin{equation}
MT_{a,b}\left(  \beta\right)  =\sup_{\alpha\in\left(  0,\alpha_{N}\right)
}\left(  \frac{1-\left(  \frac{\alpha}{\alpha_{N}}\right)  ^{\frac{N-1}{N}a}%
}{\left(  \frac{\alpha}{\alpha_{N}}\right)  ^{\frac{N-1}{N}b}}\right)
^{\frac{N-\beta}{b}}AT\left(  \alpha,\beta\right)  . \label{1.3.2}%
\end{equation}
In particular, $MT\left(  \beta\right)  <\infty$ and
\[
MT\left(  \beta\right)  =\sup_{\alpha\in\left(  0,\alpha_{N}\right)  }\left(
\frac{1-\left(  \frac{\alpha}{\alpha_{N}}\right)  ^{N-1}}{\left(  \frac
{\alpha}{\alpha_{N}}\right)  ^{N-1}}\right)  ^{\frac{N-\beta}{N}}AT\left(
\alpha,\beta\right)  .
\]

Concerning the weighted versions of the Moser-Trudinger inequalities, Ishiwata, Nakamura and Wadade \cite{INW} investigate for the scaling invariant form for the weighted Moser-Trudinger inequality by finding its best constant and proving the existence of a maximizer for the associated variational problem. Indeed, they proved the following inequality:\\
\\
\textbf{Theorem F.}  Assume $N\ge 2,-\infty<s\le t<N\quad and\quad 0<\alpha<\alpha_{N,t}:=(N-t)\omega^{\frac{1}{N-1}}_{N-1}$, then there exists a positive constant $C=C(N,s,t,\alpha)$ such that the inequality
\begin{displaymath}
\int_{\mathbb{R}^N}\Phi_N(\alpha|u|^{N/(N-1)})\frac{dx}{|x|^t}\le C(\int_{\mathbb{R}^N}\frac{|u|^{N}}{|x|^s}dx)^{\frac{N-t}{N-s}}
\end{displaymath}
holds for all radially symmetric functions $u\in L^N(\mathbb{R}^N;|x|^{-s}dx)\cap \dot{H}^{1,N}(\mathbb{R}^N)$ with $\|\nabla u\|_{L^N(\mathbb{R}^N)}\le 1$. Also the constant $\alpha_{N,t}=(N-t)\omega^{\frac{1}{N-1}}_{N-1}$ is sharp for the inequality.\\

To prove Theorem F in \cite{INW}, by taking advantage of the spherical symmetry of the functions under consideration,  they define a function $v(x):=\Big(\frac{N-t}{N}\Big)^{\frac{N-1}{N}}\tilde{u}(|x|^{\frac{N}{N-t}})$ where $u(x)=\tilde{u}(|x|)$, direct computations show that
\begin{align}
&\|\nabla u\|_{L^N(\mathbb{R}^N)}=\|\nabla v\|_{L^N(\mathbb{R}^N)} \notag\\
&\|u\|_{L^N(\mathbb{R}^N;|x|^{-t}dx)}=\frac{N}{N-t}\|v\|_{L^N(\mathbb{R}^N)} \notag\\
&\int_{\mathbb{R}^N}\Phi_N(\alpha|u|^{N/(N-1)})\frac{dx}{|x|^t}=\frac{N}{N-t}\int_{\mathbb{R}^N}\Phi_N(\frac{N}{N-t}\alpha|v|^{N/(N-1)})dx. \notag
\end{align}
Thus,  the weighted parts could be eliminated, then the proof of Theorem F  can be reduced  to that of Theorem B.  \\

The above argument cannot work if the function $u$ is not radially symmetric.
Then a natural question to ask is: can we  remove the radially symmetric condition for functions $u$ under consideration in Theorem F? We will prove in this paper that Theorem F is indeed true even when  $u$ is not necessarily radially symmetric. This is the first main result of our paper. \\

It is interesting to note that in Theorem F we could not apply the symmetrization method given by Moser in \cite{Mo} because of the existence of the weights.   \\

Now we shall state our main result of this paper, we assume the condition of exponents as follows:
\begin{equation}\label{equ1.1}
N\ge 2,-\infty<s\le t<N,N'=\frac{N}{N-1}\text{ and }0<\alpha<\alpha_{N,t}:=(N-t)\omega^{\frac{1}{N-1}}_{N-1}
\end{equation}
\\
\textbf{Theorem 1.1} \emph{Assume (\ref{equ1.1}), then there exists a positive constant $C=C(N,s,t,\alpha)$ such that the inequality
\begin{equation}\label{equ1.2}
\int_{\mathbb{R}^N}\Phi_N(\alpha|u|^{N'})\frac{dx}{|x|^t}\le C\left(\int_{\mathbb{R}^N}\frac{|u|^{N}}{|x|^s}dx\right)^{\frac{N-t}{N-s}}
\end{equation}
holds for all functions $u\in L^N(\mathbb{R}^N;|x|^{-s}dx)\cap \dot{H}^{1,N}(\mathbb{R}^N)$ with $\|\nabla u\|_{L^N(\mathbb{R}^N)}\le 1$. Moreover,  the constant $\alpha_{N,t}=(N-t)\omega^{\frac{1}{N-1}}_{N-1}$ is sharp in the sense that if $\alpha>\alpha_{N,t}$ then   the inequality (\ref{equ1.2}) cannot hold with a uniform $C$ independent of $u$.}\\

To prove  Theorem 1.1 for functions which are not necessarily radially symmetric,   we employ a different method than that of  Ishiwata, Nakamura and Wadade in \cite{INW} to prove Theorem F. The main idea is to apply a new way of change of variables to eliminate the weights in inequality (\ref{equ1.2}). We will define a new function $v$ corresponding  to $u$ which could keep the gradient norm less than 1, and eliminate the weights of integral at the same time.\\

Next,  we notice that $\Phi_N(\alpha|u|^{N/(N-1)})\le Ce^{\alpha|u|^{N/(N-1)}}|u|^N$.  Then we like to know:    could we extend the inequality in Theorem F by replacing the function $\Phi_N$ by $e^{\alpha|u|^{N/(N-1)}}|u|^N$ on the left hand side?
This is the second main result of this paper.\\
\\
\textbf{Theorem 1.2} \emph{Assume (\ref{equ1.1}), then there exists a positive constant $C=C(N,s,t,\alpha)$ such that the inequality
\begin{equation}\label{equ1.3}
\int_{\mathbb{R}^N}e^{\alpha|u|^{N'}}|u|^N\frac{dx}{|x|^t}\le C\left(\int_{\mathbb{R}^N}\frac{|u|^{N}}{|x|^s}dx\right)^{\frac{N-t}{N-s}}
\end{equation}
holds for all functions $u\in L^N(\mathbb{R}^N;|x|^{-s}dx)\cap \dot{H}^{1,N}(\mathbb{R}^N)$ with $\|\nabla u\|_{L^N(\mathbb{R}^N)}\le 1$. Moreover,  the constant $\alpha_{N,t}=(N-t)\omega^{\frac{1}{N-1}}_{N-1}$ is sharp in the sense that if $\alpha>\alpha_{N,t}$ then the inequality (\ref{equ1.3}) cannot hold with a uniform $C$ independent of $u$.}\\

To prove  Theorem 1.2, we verify the non-singular case, which states that, for  $N\ge2$ and $0<\alpha<\alpha_N$, there exists a positive constant $C=C(N,\alpha)$ such that the inequality
\begin{displaymath}
\int_{\mathbb{R}^N}e^{\alpha|u|^{N'}}|u|^Ndx\le C\|u\|^N_{L^N(\mathbb{R}^N)}
\end{displaymath}
holds for all $u\in W^{1,N}(\mathbb{R}^N)$ with $\|\nabla u\|_{L^N(\mathbb{R}^N)}\le1$. Then (\ref{equ1.3}) can be obtained by using the same method of changing variables  used in the proof of  Theorem 1.1 to eliminate the weights.\\

In fact, for $q>N$, we could have a more general form for this inequality.\\
\\
\textbf{Theorem 1.3} \emph{Assume (\ref{equ1.1}), then there exists a positive constant $C=C(N,q,s,t,\alpha)$ such that the inequality
\begin{equation}\label{equ1.4}
\int_{\mathbb{R}^N}e^{\alpha|u|^{N'}}|u|^q\frac{dx}{|x|^t}\le\Big(\int_{\mathbb{R}^N}\frac{|u|^q}{|x|^s}dx\Big)^{\frac{N-t}{N-s}}
\end{equation}
holds for all functions $u\in L^q(\mathbb{R}^N;|x|^{-s}ds)\cap \dot{H}^{1,N}(\mathbb{R})$ with $\|\nabla u\|_{L^N(\mathbb{R}^N)}\le 1$.}\\

Next, we shall discuss the existence of a maximizer associated with each of our inequalities. Ishiwata, Nakamura and Wadade have proved in \cite{INW} the existence of an maximizer for the inequality (\ref{equ1.2}) in Theorem 1.1   for radially symmetric functions.\\

To extend this result to functions that are not necessarily symmetric, we will show that any maximizing sequence must be obtained when they are radially symmetric, consequently, we only need to consider the radially symmetric functions.\\

Use $X^{1,N}_s$ and $X^{1,N}_{s,rad}$ denote the weighted Sobolev spaces defined by
\begin{align}
\left\{
\begin{array}{rcl}
X^{1,N}_s&:=&L^N(\mathbb{R}^N;|x|^{-s}ds)\cap \dot{H}^{1,N}(\mathbb{R}), \\
X^{1,N}_{s,rad}&:=&\{u\in X^{1,N}_s, u \text{ is radially symmetric}\}. \notag
\end{array}
\right.
\end{align}

Then we define the sharp constants $\mu_{N,s,t,\alpha}(\mathbb{R}^N)$ and $\nu_{N,s,t,\alpha}(\mathbb{R}^N)$ of (\ref{equ1.2}) and (\ref{equ1.3}) by
\begin{displaymath}
\mu_{N,s,t,\alpha}(\mathbb{R}^N):=\sup_{\substack{u\in X^{1,N}_s\\\|\nabla u\|_{L^N(\mathbb{R}^N)}=1}} F_{N,s,t,\alpha}(u),
\end{displaymath}
\begin{displaymath}
\nu_{N,s,t,\alpha}(\mathbb{R}^N):=\sup_{\substack{u\in X^{1,N}_s\\\|\nabla u\|_{L^N(\mathbb{R}^N)}=1}} G_{N,s,t,\alpha}(u),
\end{displaymath}
where
\begin{equation}\label{equ1.5}
F_{N,s,t,\alpha}(u):=\frac{\int_{\mathbb{R}^N}\Phi_N(\alpha|u|^{N'})\frac{dx}{|x|^t}}{\|u\|^{\frac{N(N-t)}{N-s}}_{L^N(\mathbb{R}^N;|x|^{-s}dx)}},
\end{equation}
\begin{equation}\label{equ1.6}
G_{N,s,t,\alpha}(u):=\frac{\int_{\mathbb{R}^N}e^{\alpha|u|^{N'}}|u|^N\frac{dx}{|x|^t}}{\|u\|^{\frac{N(N-t)}{N-s}}_{L^N(\mathbb{R}^N;|x|^{-s}dx)}}.
\end{equation}

By a suitable renormalization argument and compact embedding theorem for radial Sobolev space we prove the following Theorem.\\
\\
\textbf{Theorem 1.4} (i)\emph{Assume (\ref{equ1.1}) holds, then the sharp constant $\mu_{N,s,t,\alpha}(\mathbb{R}^N)$ is attained.}

(ii) \emph{Assume (\ref{equ1.1}) holds, then the sharp constant $\nu_{N,s,t,\alpha}(\mathbb{R}^N)$ is attained.}\\

The paper is organized as follows: In Section 2, to eliminate the weights in the weighted Moser-Trudinger inequality in Theorem 1.1, we will employ a new method of change of variables   to establish Theorem 1.1. In Section 3 we prove two lemmas directly  corresponding  to Theorem 1.2 and Theorem 1.3. Then we will complete the proof of Theorem 1.2 and Theorem 1.3 in Section 4. The existence of the maximizer (Theorem 1.4)  will be established in Section 5.\\

\section{Proof of Theorem 1.1}

 It is not hard to see that it suffices to prove that  inequality (\ref{equ1.2}) holds for the special case $s=t$, which states that, under the assumption (\ref{equ1.1}), there exists a positive constant $C=C(N,t,\alpha)$ such that the inequality
\begin{equation}\label{equ2.1}
\int_{\mathbb{R}^N}\Phi_N(\alpha|u|^{N'})\frac{dx}{|x|^t}\le C\int_{\mathbb{R}^N}\frac{|u|^N}{|x|^t}dx
\end{equation}
holds for all functions $u\in X^{1,N}_s$ with $\|\nabla u\|_{L^N(\mathbb{R}^N)}\le 1$.\\

Once we have proved this special case (\ref{equ2.1}), the general case $s<t$ follows immediately by applying the following Caffarelli-Kohn-Nirenberg inequality established in \cite{CKN}.  \\

\textbf{Theorem G. (Caffarelli-Kohn-Nirenberg inequality)} \emph{For $u\in C^{\infty}_{0}(\mathbb{R}^n)$. In what follows $p,q,r;\alpha,\beta,\sigma$ and $a$ are fixed real numbers satisfying
\begin{equation}\label{equ2.2}
p,q\ge 1,\quad r>0,\quad 0\le a\le 1,
\end{equation}
\begin{equation}\label{equ2.3}
\frac{1}{p}+\frac{\alpha}{n},\quad \frac{1}{q}+\frac{\beta}{n},\quad \frac{1}{r}+\frac{\gamma}{n}>0,
\end{equation}
where $\gamma=a\sigma+(1-a)\beta$.}

\emph{There exists a positive constant $C$ such that the following inequality holds for all $u\in C^{\infty}_{0}(\mathbb{R}^n)$,
\begin{equation}\label{equ2.4}
\big{\|}|x|^{\gamma}u\big{\|}_{L^r}\le C\big{\|}|x|^{\alpha}|Du|\big{\|}^{a}_{L^p}\big{\|}|x|^{\beta}u\big{\|}^{1-a}_{L^q},
\end{equation}
if and only if the following relations hold:
\begin{equation}\label{equ2.5}
\frac{1}{r}+\frac{\gamma}{n}=a(\frac{1}{p}+\frac{\alpha-1}{n})+(1-a)(\frac{1}{q}+\frac{\beta}{n}).
\end{equation}}

\emph{Furthermore, on any compact set in the parameter space in which (\ref{equ2.2}), (\ref{equ2.3}), (\ref{equ2.5}) and $0\le\alpha-\sigma\le 1$ hold, the constant $F$ is bounded.}\\

Therefore for $q\ge N$, applying the conditions in this Theorem we have
\begin{equation}\label{equ2.6}
\|u\|_{L^q(\mathbb{R}^N;|x|^{-t}dx)}\le C\|u\|^{\frac{N-t}{N-s}}_{L^q(\mathbb{R}^N;|x|^{-s}dx)}\|\nabla u\|^{1-\frac{N-t}{N-s}}_{L^N(\mathbb{R}^N)}.
\end{equation}

Apply this to (\ref{equ2.1}) we can directly get the inequality (\ref{equ1.2}) in Theorem 1.1.\\

Now we begin the proof of Theorem 1.1.

\begin{proof}
Let $0<\alpha<\alpha_{N,t}$ and let $u\in X^{1,N}_s$ with $\|\nabla u\|_{L^N(\mathbb{R}^N)}\le 1$. We define the function $v\in W^{1,N}(\mathbb{R}^N)$ for $x\in\mathbb{R}^N$ by the formula below,
\begin{equation}\label{equ2.7}
v(x):=\Big(\frac{N-t}{N}\Big)^{\frac{1}{N'}}u(|x|^{\frac{t}{N-t}}x).
\end{equation}

Consider the vector-valued function $F:\mathbb{R}^N\to\mathbb{R}^N$ defined by
\begin{displaymath}
F(x)=|x|^{\frac{t}{N-t}}x,
\end{displaymath}
the Jacobian matrix of this function $F$ is
\begin{displaymath}
\mathbf{J_F} =
\left( \begin{array}{cccc}
|x|^{\frac{t}{N-t}}+\frac{t}{N-t}x^2_1|x|^{\frac{3t-2N}{N-t}} & \frac{t}{N-t}x_1 x_2|x|^{\frac{3t-2N}{N-t}} & \ldots & \frac{t}{N-t}x_1 x_N|x|^{\frac{3t-2N}{N-t}} \\
\frac{t}{N-t}x_2 x_1|x|^{\frac{3t-2N}{N-t}} & |x|^{\frac{t}{N-t}}+\frac{t}{N-t}x^2_2|x|^{\frac{3t-2N}{N-t}} & \ldots & \frac{t}{N-t}x_2 x_N|x|^{\frac{3t-2N}{N-t}} \\
\vdots & \vdots & \ddots & \vdots \\
\frac{t}{N-t}x_N x_1|x|^{\frac{3t-2N}{N-t}} & \frac{t}{N-t}x_N x_2|x|^{\frac{3t-2N}{N-t}} & \ldots & |x|^{\frac{t}{N-t}}+\frac{t}{N-t}x^2_N|x|^{\frac{3t-2N}{N-t}}
\end{array} \right),
\end{displaymath}
direct calculations show us
\begin{equation}\label{equ2.8}
\det(J_F)=\frac{N}{N-t}|x|^{\frac{Nt}{N-t}}.
\end{equation}

Then for
\begin{displaymath}
\int_{\mathbb{R}^N}|v(x)|^N dx=\Big(\frac{N-t}{N}\Big)^{N-1}\int_{\mathbb{R}^N}\big|u(|x|^{\frac{t}{N-t}}x)\big|^Ndx,
\end{displaymath}
using change of variables $y_i=|x|^{\frac{t}{N-t}}x_i$, $i=1,2,...,N$. We have
\begin{equation}\label{equ2.9}
dy=\det(J_F)dx=\frac{N}{N-t}|x|^{\frac{Nt}{N-t}}dx,
\end{equation}
and
\begin{equation}\label{equ2.10}
dx=\frac{N-t}{N}\frac{dy}{|y|^t},
\end{equation}
therefore,   we have
\begin{align}\label{equ2.11}
\int_{\mathbb{R}^N}|v(x)|^N dx&=\Big(\frac{N-t}{N}\Big)^{N-1}\int_{\mathbb{R}^N}\big|u(|x|^{\frac{t}{N-t}}x)\big|^N dx \notag\\
&=\Big(\frac{N-t}{N}\Big)^N\int_{\mathbb{R}^N}|u(y)|^N\frac{dy}{|y|^t}.
\end{align}

Now we begin to consider the gradient of $v$. After calculations, we have
\begin{align}
\left( \begin{array}{cccc}
\frac{\partial v}{\partial x_1}(x) \\
\frac{\partial v}{\partial x_2}(x) \\
\vdots \\
\frac{\partial v}{\partial x_N}(x)
\end{array} \right)
=\nabla v(x)=\Big(\frac{N-t}{N}\Big)^{\frac{1}{N'}}\nabla(u(|x|^{\frac{t}{N-t}}x))=\Big(\frac{N-t}{N}\Big)^{\frac{1}{N'}}
J^T_F\left( \begin{array}{cccc}
\frac{\partial u}{\partial x_1}(|u|^{\frac{t}{N-t}}x) \\
\frac{\partial u}{\partial x_2}(|u|^{\frac{t}{N-t}}x) \\
\vdots \\
\frac{\partial u}{\partial x_N}(|u|^{\frac{t}{N-t}}x)
\end{array} \right), \notag
\end{align}
Hence we have
\begin{displaymath}
\frac{\partial v}{\partial x_i}(x)=\Big(\frac{N-t}{N}\Big)^{\frac{1}{N'}}\Big(|x|^{\frac{t}{N-t}}\frac{\partial u}{\partial x_i}(|x|^{\frac{t}{N-t}}x)+A_i\Big),
\end{displaymath}
for $i=1,2,...N$, where $A_i$ is defined by following
\begin{displaymath}
A_i:=\sum^N_{j=1}\frac{t}{N-t}x_i x_j|x|^{\frac{3t-2N}{N-t}}\frac{\partial u}{\partial x_j}(|x|^{\frac{t}{N-t}}x),
\end{displaymath}
 Substituting  them into $|\nabla v(x)|^2$, we obtain
\begin{align}
|\nabla v(x)|^2&= \sum^N_{i=1}\Big(\frac{\partial v}{\partial x_i}(x)\Big)^2 \notag\\
&=\Big(\frac{N-t}{N}\Big)^{\frac{2}{N'}}\sum^N_{i=1}\Big(|x|^{\frac{t}{N-t}}\frac{\partial u}{\partial x_i}(|x|^{\frac{t}{N-t}}x)+A_i\Big)^2 \notag\\
&=\Big(\frac{N-t}{N}\Big)^{\frac{2}{N'}}\Big(\sum^N_{i=1}|x|^{\frac{2t}{N-t}}\Big(\frac{\partial u}{\partial x_i}(|x|^{\frac{t}{N-t}}x)\Big)^2
+\sum^N_{i=1}2A_i|x|^{\frac{t}{N-t}}\frac{\partial u}{\partial x_i}(|x|^{\frac{t}{N-t}}x)+\sum^N_{i=1}A^2_i\Big) \notag\\
&:=\Big(\frac{N-t}{N}\Big)^{\frac{2}{N'}}\Big(I_1+I_2+I_3\Big). \notag
\end{align}

Direct computations show us the first term
\begin{align}
I_1=|x|^{\frac{2t}{N-t}}\big|\nabla u(|x|^{\frac{t}{N-t}}x)\big|^2, \notag
\end{align}
Applying the Cauchy-Schwarz inequality to estimate the second term, we  get
\begin{align}
I_2&=\sum^N_{i=1}2A_i|x|^{\frac{t}{N-t}}\frac{\partial u}{\partial x_i}(|x|^{\frac{t}{N-t}}x) \notag\\
&=\frac{2t}{N-t}|x|^{\frac{2t}{N-t}}\sum^N_{i=1}\sum^N_{j=1}\frac{x_i x_j}{|x|^2}\frac{\partial u}{\partial x_j}(|x|^{\frac{t}{N-t}}x)
\frac{\partial u}{\partial x_i}(|x|^{\frac{t}{N-t}}x) \notag\\
&=\frac{2t}{N-t}|x|^{\frac{2t}{N-t}}\Big(\sum^N_{i=1}\frac{x_i}{|x|}\frac{\partial u}{\partial x_i}(|x|^{\frac{t}{N-t}}x)\Big)^2 \notag\\
&\le\frac{2t}{N-t}|x|^{\frac{2t}{N-t}}\Big(\sum^N_{i=1}\big(\frac{x_i}{|x|}\big)^2\Big)\Big(\sum^N_{i=1}\big(\frac{\partial u}{\partial x_i}(|x|^{\frac{t}{N-t}}x)\big)^2\Big) \notag\\
&=\frac{2t}{N-t}|x|^{\frac{2t}{N-t}}\big|\nabla u(|x|^{\frac{t}{N-t}}x)\big|^2, \notag
\end{align}
Similarly for the last term we have
\begin{align}
I_3=\sum^N_{i=1}A^2_i&=\sum^N_{i=1}\Big(\sum^N_{j=1}\frac{t}{N-t}x_i x_j|x|^{\frac{3t-2N}{N-t}}\frac{\partial u}{\partial x_j}(|x|^{\frac{t}{N-t}}x)\Big)^2 \notag\\
&\le\Big(\frac{t}{N-t}\Big)^2\sum^N_{i=1}\Bigg(\Big(\sum^N_{j=1}\big(x_i x_j|x|^{\frac{3t-2N}{N-t}}\big)^2\Big)
\Big(\sum^N_{j=1}\big(\frac{\partial u}{\partial x_j}(|x|^{\frac{t}{N-t}}x)\big)^2\Big)\Bigg) \notag\\
&=\Big(\frac{t}{N-t}\Big)^2|x|^{\frac{2t}{N-t}}\big|\nabla u(|x|^{\frac{t}{N-t}}x)\big|^2, \notag
\end{align}
Combining them together we have
\begin{displaymath}
|\nabla v(x)|^2\le\Big(\frac{N-t}{N}\Big)^{\frac{2}{N'}}\Big(|x|^{\frac{2t}{N-t}}+\frac{2t}{N-t}|x|^{\frac{2t}{N-t}}+\big(\frac{t}{N-t}\big)^2|x|^{\frac{2t}{N-t}}\Big)\big|\nabla u(|x|^{\frac{t}{N-t}}x)\big|^2.
\end{displaymath}
This  leads to
\begin{displaymath}
|\nabla v(x)|\le\Big(\frac{N}{N-t}\Big)^{\frac{1}{N}}|x|^{\frac{t}{N-t}}\big|\nabla u(|x|^{\frac{t}{N-t}}x)\big|.
\end{displaymath}

Using the change of variables again, we get
\begin{align}\label{equ2.12}
\int_{\mathbb{R}^N}|\nabla v(x)|^N dx&\le\frac{N}{N-t}\int_{\mathbb{R}^N}|x|^{\frac{Nt}{N-t}}\big|\nabla u(|x|^{\frac{t}{N-t}}x)\big|^N dx \notag\\
&=\int_{\mathbb{R}^N}|\nabla u(y)|^N dy.
\end{align}

By (\ref{equ2.11}),(\ref{equ2.12}) we obtain
\begin{equation}\label{equ2.13}
\|u\|_{L^N(\mathbb{R}^N;|x|^{-t}dx)}=\frac{N}{N-t}\|v\|_{L^N(\mathbb{R}^N)},
\end{equation}
\begin{equation}\label{equ2.14}
\|\nabla u\|_{L^N(\mathbb{R}^N)}\ge\|\nabla v\|_{L^N(\mathbb{R}^N)}.
\end{equation}

From computations we also have
\begin{align}\label{equ2.15}
&\int_{\mathbb{R}^N}\Phi_N(\alpha|u(y)|^{N'})\frac{dy}{|y|^t} \notag\\
=&\sum^{\infty}_{i=N-1}\int_{\mathbb{R}^N}\frac{(\alpha|u(y)|^{N'})^i}{i!}\frac{dy}{|y|^t} \notag\\
=&\sum^{\infty}_{i=N-1}\int_{\mathbb{R}^N}\frac{N}{N-t}\frac{\alpha^i}{i!}\big|u(|x|^{\frac{t}{N-t}}x)\big|^{N'i}dx \notag\\
=&\sum^{\infty}_{i=N-1}\int_{\mathbb{R}^N}\Big(\frac{N}{N-t}\Big)^{i+1}\frac{\alpha^i}{i!}|v(x)|^{N'i}dx \notag\\
=&\frac{N}{N-t}\sum^{\infty}_{i=N-1}\int_{\mathbb{R}^N}\frac{(\frac{N}{N-t}\alpha|v(x)|^{N'})^i}{i!}dx \notag\\
=&\frac{N}{N-t}\int_{\mathbb{R}^N}\Phi_N\big(\frac{N}{N-t}\alpha|v(x)|^{N'}\big)dx.
\end{align}

Since $0<\frac{N}{N-t}\alpha<\frac{N}{N-t}\alpha_{N,t}=N\omega^{\frac{1}{N-1}}_{N-1}$ and $\|\nabla v\|_{L^N(\mathbb{R}^N)}\le 1$, by applying Theorem B and (\ref{equ2.13}), (\ref{equ2.15}) we have
\begin{align}\label{equ2.16}
\int_{\mathbb{R}^N}\Phi_N(\alpha|u(x)|^{N'})\frac{dx}{|x|^t} \notag=&\frac{N}{N-t}\int_{\mathbb{R}^N}\Phi_N\big(\frac{N}{N-t}\alpha|v(x)|^{N'}\big)dx \notag\\
\le&C_{\alpha}\|v\|^N_{L^N(\mathbb{R}^N)} \notag\\
=&C_{\alpha}\|u\|^N_{L^N(\mathbb{R}^N;|x|^{-t}dx)}.
\end{align}
This is exactly the special case (\ref{equ2.1}), therefore we have proved the inequality (\ref{equ1.2}) in Theorem 1.1 by using the Caffarelli-Kohn-Nirenberg inequality.\\

Next,  we want to show $\alpha_{N,t}$ is the sharp constant for (\ref{equ1.2}), here we apply the following modified Moser's test sequence  used in \cite{INW}.\\

For $k\in\mathbb{N}$, define a sequence $u_k$ such that $u_k\in X^{1,N}_s$ by
\begin{align}\label{equ2.17}
u_k(x)=\left\{
\begin{aligned}
&0 &if\quad |x|\ge 1 \\
&(\frac{N-t}{\omega_{N-1}k})^{\frac{1}{N}}\log{\frac{1}{|x|}}& \text{ if } e^{-\frac{k}{N-t}}<|x|<1 \\
&(\frac{1}{\omega_{N-1}})^{\frac{1}{N}}(\frac{k}{N-t})^{\frac{1}{N'}}& \text{ if } 0\le|x|\le e^{-\frac{k}{N-t}}
\end{aligned}
\right..
\end{align}

Direct computation show that $\|\nabla u_k\|_{L^N(\mathbb{R}^N)}=1$ for all $k\in\mathbb{N}$, and we have
\begin{align}
\frac{\int_{\mathbb{R}^N}e^{\alpha_{N,t}|u_k|^{N'}}|u_k|^N\frac{dx}{|x|^t}}{\|u_k\|^{\frac{N(N-t)}{N-s}}_{L^N(\mathbb{R}^N;|x|^{-s}dx)}}\to\infty \text{   as }k\to\infty, \notag
\end{align}
which implies inequality (\ref{equ1.2}) fails when $\alpha=\alpha_{N,t}$, hence we finish the proof of Theorem 1.1.
\end{proof}

\section{Two lemmas}

In this section, we provide the proof for two lemmas in non-singular form (e.g. inequality (1.3) when $t=s=0$). The proofs of these two lemmas can be done using
an idea used in the work of Lam and the second author \cite{LamLu} by considering level sets of the functions under consideration.
This can be carried out   in more general singular case (including the case $s=0$, but $t\not =0$) without using symmetrization. However, we present a proof using the symmetrization argument of Moser \cite{Mo} in the non-singular $t=s=0$. \\
\\
\textbf{Lemma 3.1} \emph{Suppose $N\ge 2$. Then for any $\alpha\in(0,\alpha_{N})$, where $\alpha_N=N\omega^{\frac{1}{N-1}}_{N-1}$, there exists a constant $C_{\alpha}>0$ such that
\begin{equation}\label{equ3.1}
\int_{\mathbb{R}^N}e^{\alpha|u|^{N'}}|u|^N dx\le C_{\alpha}\|u\|^N_{L^N(\mathbb{R}^N)}
\end{equation}
holds for all functions $u\in W^{1,N}(\mathbb{R}^N)$ with $\|\nabla u\|_{L^N(\mathbb{R}^N)}\le 1$.}

\begin{proof}
To prove this lemma, we use the idea of means of symmetrization given by Moser \cite{Mo}. Then it is suffices for us to show inequality (\ref{equ3.1}) satisfied for non-negative, compactly supported, radially symmetric functions $u(x)=\tilde{u}(|x|)$, and $\tilde{u}(|x|):[0,\infty)\to\mathbb{R}$ are decreasing.\\

Following Moser's argument, we set
\begin{displaymath}
w(t)=N^{\frac{N-1}{N}}\omega^{\frac{1}{N}}_{N-1}\tilde{u}(e^{-\frac{t}{N}}),\quad\quad |x|^N=e^{-t}.
\end{displaymath}
Then we have $w(t)$ defined on $(-\infty,\infty)$ and satisfied
\begin{equation}\label{equ3.2}
w(t)\ge 0\quad\text{ for }t\in\mathbb{R},
\end{equation}
\begin{equation}\label{equ3.3}
w'(t)\ge 0\quad\text{ for }t\in\mathbb{R},
\end{equation}
\begin{equation}\label{equ3.4}
w(t_0)=0\quad\text{ for some }t_0\in\mathbb{R}.
\end{equation}
From calculation we have
\begin{align}\label{equ3.5}
\int_{\mathbb{R}^N}|\nabla u|^N dx&=\int^{\infty}_0\omega_{N-1}|\tilde{u}'(|x|)|^N |x|^{N-1}d|x| \notag\\
&=\int^{\infty}_{-\infty}N^{N-1}\omega_{N-1}|\tilde{u}'(e^{-\frac{t}{N}})|\frac{e^{-t}}{N^N}dt \notag\\
&=\int^{\infty}_{-\infty}|w'(t)|^N dt.
\end{align}
\begin{align}\label{equ3.6}
\int_{\mathbb{R}^N}|u(x)|^N dx&=\int^{\infty}_0\omega_{N-1}|\tilde{u}(|x|)|^N |x|^{N-1}d|x| \notag\\
&=\frac{1}{N^N}\int^{\infty}_{-\infty}N^{N-1}\omega_{N-1}|\tilde{u}(e^{-\frac{t}{N}})|^N e^{-t}dt \notag\\
&=\frac{1}{N^N}\int^{\infty}_{-\infty}|w(t)|^N e^{-t} dt.
\end{align}
\begin{align}\label{equ3.7}
\int_{\mathbb{R}^N}e^{\alpha|u(x)|^{N'}}|u(x)|^N dx=&\sum^{\infty}_{j=N-1}\frac{\alpha^{j-(N-1)}}{(j-(N-1))!}\int_{\mathbb{R}^N}(|u(x)|^{N'})^j dx \notag\\
=&\sum^{\infty}_{j=N-1}\frac{\alpha^{j-(N-1)}}{(j-(N-1))!}\int^{\infty}_0 \omega_{N-1}(|\tilde{u}(|x|)|^{N'})^j |x|^{N-1}d|x| \notag\\
=&\frac{\omega_{N-1}}{N}\sum^{\infty}_{j=N-1}\frac{\alpha^{j-(N-1)}}{(j-(N-1))!}\int^{\infty}_{-\infty}\big(\frac{1}{\alpha_N}|w(t)|^{N'}\big)^j e^{-t}dt \notag\\
=&\frac{1}{N^N}\int^{\infty}_{-\infty}e^{\frac{\alpha}{\alpha_N}|w(t)|^{N'}}|w(t)|^N e^{-t}dt.
\end{align}

Therefore to prove our lemma it suffices to show that for $\beta\in(0,1)$ there exists $C_{\beta}>0$ such that
\begin{equation}\label{equ3.8}
\int^{\infty}_{-\infty}e^{\beta|w(t)|^{N'}}|w(t)|^N e^{-t}dt\le C_{\beta}\int^{\infty}_{-\infty}|w(t)|^N e^{-t}dt
\end{equation}
for all function $w(t)$ satisfying the conditions (\ref{equ3.2})-(\ref{equ3.4}) and $\int^{\infty}_{-\infty}|w'(t)|^N dt=1$.\\

Set $T_0=sup\{t\in\mathbb{R}|w(t)\le 1\}\in(-\infty,\infty]$, then we split the integral set to be $(-\infty,T_0]\cup[T_0,\infty)$. Next we will show the inequality satisfied for each of them.\\

For $t\in(-\infty,T_0]$, we have $w(t)\in[0,1]$, therefore $e^{\beta|w(t)|^{N'}}\le e^{\beta}:=C_{1,\beta}$ on this integral part. Hence we have
\begin{equation}\label{equ3.9}
\int^{T_0}_{-\infty}e^{\beta|w(t)|^{N'}}|w(t)|^N e^{-t}dt\le C_{1,\beta}\int^{T_0}_{-\infty}|w(t)|^N e^{-t}dt.
\end{equation}

Then we consider the second integral over $[T_0,\infty)$. Since $w(T_0)=1$, apply H\"{o}lder's inequality we have for $t\ge T_0$
\begin{align}
w(t)&=w(T_0)+\int^t_{T_0}w'(\tau)d\tau \notag\\
&\le 1+(t-T_0)^{\frac{1}{N'}}\Big(\int^{\infty}_{T_0}w'(\tau)^N d\tau\Big)^{\frac{1}{N}} \notag\\
&\le 1+(t-T_0)^{\frac{1}{N'}}. \notag
\end{align}

Then we need to apply an inequality, for any $\epsilon>0$ there exists a constant $C_{\epsilon}>0$ such that
\begin{align}
1+s^{\frac{1}{N'}}\le((1+\epsilon)s+C_{\epsilon})^{\frac{1}{N'}} \notag
\end{align}
for all $s\ge 0$. Therefore we have
\begin{equation}\label{equ3.10}
w(t)\le [(1+\epsilon)(t-T_0)+C_{\epsilon}]^{\frac{1}{N'}}.
\end{equation}

For any $\beta\in(0,1)$, it is possible for us to choose an $\epsilon$ small enough such that $\beta(1+\epsilon)<1$.

Then applying (\ref{equ3.10}) for our integral we have
\begin{align}
&\int^{\infty}_{T_0}e^{\beta|w(t)|^{N'}}|w(t)|^N e^{-t}dt \notag\\
\le&\int^{\infty}_{T_0}\exp\big(\beta((1+\epsilon)(t-T_0)+C_{\epsilon})-t\big)\big((1+\epsilon)(t-T_0)+C_{\epsilon}\big)^{N-1}dt \notag\\
=&\int^{\infty}_{T_0}\exp\big((\beta(1+\epsilon)-1)(t-T_0)+\beta C_{\epsilon}-T_0\big)\big((1+\epsilon)(t-T_0)+C_{\epsilon}\big)^{N-1}dt \notag\\
=&I_1+J_1, \notag
\end{align}
where $I_1$ and $J_1$ are obtained by using integration by parts as follows.\\

For short we set $A_{\beta}(t)=\exp\big((\beta(1+\epsilon)-1)(t-T_0)+\beta C_{\epsilon}-T_0\big)$, therefore we get
\begin{align}
I_1=\frac{1}{\beta(1+\epsilon)-1}A_{\beta}(t)\big((1+\epsilon)(t-T_0)+C_{\epsilon}\big)^{N-1}\Big|^{t=\infty}_{t=T_0}
=\frac{C^{N-1}_{\epsilon}e^{\beta C_{\epsilon}}}{1-\beta(1+\epsilon)}e^{-T_0}, \notag
\end{align}
\begin{align}
J_1=\frac{(N-1)(1+\epsilon)}{1-\beta(1+\epsilon)}\int^{\infty}_{T_0}A_{\beta}(t)\big((1+\epsilon)(t-T_0)+C_{\epsilon}\big)^{N-2}dt. \notag
\end{align}

We could apply the similar integration by parts repeatedly and define
\begin{align}
&J_k=I_{k+1}+J_{k+1},\quad\quad k=1,2,...,N-1. \notag
\end{align}
Hence we have
\begin{displaymath}
I_1+J_1=I_1+I_2+J_2=...=J_N+\sum^N_{k=1}I_k,
\end{displaymath}
direct computation shows us
\begin{equation}\label{equ3.11}
J_k=\frac{(N-1)!(1+\epsilon)^k}{(N-k-1)!(1-\beta(1+\epsilon))^k}\int^{\infty}_{T_0}A_{\beta}(t)\big((1+\epsilon)(t-T_0)+C_{\epsilon}\big)^{N-k-1}dt,
\end{equation}
for $k=1,2,...,N-1$. Therefore we have
\begin{align}
J_{N-1}=\frac{(N-1)!(1+\epsilon)^{N-1}}{(1-\beta(1+\epsilon))^{N-1}}\int^{\infty}_{T_0}A_{\beta}(t)dt, \notag
\end{align}
where it is easy to see $J_N=0$.\\

On the other hand we have
\begin{equation}\label{equ3.12}
I_k=\frac{(N-1)!(1+\epsilon)^{k-1}C^{N-k}_{\epsilon}e^{\beta C_{\epsilon}}}{(N-k)!(1-\beta(1+\epsilon))^k}e^{-T_0},
\end{equation}
for $k=1,2,...,N$.\\

Since $\epsilon$ is only related to $\beta$, we define $C_{2,\beta}$ as follows
\begin{displaymath}
C_{2,\beta}=\sum^{N}_{k=1}\frac{(N-1)!(1+\epsilon)^{k-1}C^{N-k}_{\epsilon}e^{\beta C_{\epsilon}}}{(N-k)!(1-\beta(1+\epsilon))^k}.
\end{displaymath}

Since $w(t)\ge 1$ on $[T_0,\infty)$, we have
\begin{equation}\label{equ3.13}
\int^{\infty}_{T_0}|w(t)|^N e^{-t}dt\ge\int^{\infty}_{T_0}e^{-t}dt=e^{-T_0},
\end{equation}
So we get
\begin{equation}\label{equ3.14}
\int^{\infty}_{T_0}e^{\beta|w(t)|^{N'}}|w(t)|^N e^{-t}dt\le\sum^N_{k=1}I_k=C_{2,\beta}e^{-T_0}\le C_{2,\beta}\int^{\infty}_{T_0}|w(t)|^N e^{-t}dt.
\end{equation}

Now setting $C_{\beta}=\max\{C_{1,\beta},C_{2,\beta}\}$, which is only dependent on  $\beta$ and $N$, and combining (\ref{equ3.9}) with (\ref{equ3.14}),  we get
\begin{equation}\label{equ3.15}
\int^{\infty}_{-\infty}e^{\beta|w(t)|^{N'}}|w(t)|^N e^{-t}dt\le C_{\beta}\int^{\infty}_{\infty}|w(t)|^N e^{-t}dt.
\end{equation}
Thus,  the lemma is proved.
\end{proof}

\textbf{Lemma 3.2} \emph{Assume $2\le N<q$. Then for any $\alpha\in(0,\alpha_{N})$, where $\alpha_N=N\omega^{\frac{1}{N-1}}_{N-1}$, there exists a constant $C=C(N,\alpha,q)>0$ such that
\begin{equation}\label{equ3.16}
\int_{\mathbb{R}^N}e^{\alpha|u|^{N'}}|u|^q dx\le C\|u\|^q_{L^q(\mathbb{R}^N)}
\end{equation}
holds for all functions $u\in L^q(\mathbb{R}^N)\cap H^{1,N}(\mathbb{R}^N)$ with $\|\nabla u\|_{L^N(\mathbb{R}^N)}\le 1$.}

\begin{proof}
As in Lemma 2.1, we apply the method  of symmetrization and define
\begin{displaymath}
w(t)=N^{\frac{N-1}{N}}\omega^{\frac{1}{N}}_{N-1}\tilde{u}(e^{-\frac{t}{N}}),\quad\quad |x|^N=e^{-t}.
\end{displaymath}
on $(-\infty,\infty)$ that satisfy (\ref{equ3.2})-(\ref{equ3.4}).\\

Then direct calculations show
\begin{equation}\label{equ3.17}
\int_{\mathbb{R}^N}|\nabla u|^N dx=\int^{\infty}_{-\infty}|w'(t)|^N dt,
\end{equation}
\begin{equation}\label{equ3.18}
\int_{\mathbb{R}^N}|u(x)|^q dx=N^{\frac{q}{N}-1-q}\omega^{1-\frac{q}{N}}_{N-1}\int^{\infty}_{-\infty}|w(t)|^q e^{-t} dt,
\end{equation}
\begin{equation}\label{equ3.19}
\int_{\mathbb{R}^N}e^{\alpha|u(x)|^{N'}}|u(x)|^q dx=N^{\frac{q}{N}-1-q}\omega^{1-\frac{q}{N}}_{N-1}\int^{\infty}_{-\infty}e^{\frac{\alpha}{\alpha_N}|w(t)|^{N'}}|w(t)|^q e^{-t}dt.
\end{equation}\\

Therefore, to prove our lemma it suffices to prove that for $\beta\in(0,1)$ there exists $C_{\beta}>0$ such that
\begin{equation}\label{equ3.20}
\int^{\infty}_{-\infty}e^{\beta|w(t)|^{N'}}|w(t)|^N e^{-t}dt\le C_{\beta}\int^{\infty}_{-\infty}|w(t)|^N e^{-t}dt
\end{equation}
for all function $w(t)$ satisfying the conditions (\ref{equ3.2})-(\ref{equ3.4}) and $\int^{\infty}_{-\infty}|w'(t)|^N dt=1$.\\

Arguing similarly to  what we did  in Lemma 3.1,  we set $T_0=sup\{t\in\mathbb{R}|w(t)\le 1\}\in(-\infty,\infty]$, and split the integral set $(-\infty,\infty)$ to be $(-\infty,T_0]\cup[T_0,\infty)$.\\

For $t\in(-\infty,T_0]$ we have $w(t)\in[0,1]$, therefore $e^{\beta|w(t)|^{N'}}\le e^{\beta}:=C_{1,\beta}$ on this integral part. Hence we have
\begin{equation}\label{equ3.21}
\int^{T_0}_{-\infty}e^{\beta|w(t)|^{N'}}|w(t)|^q e^{-t}dt\le C_{1,\beta}\int^{T_0}_{-\infty}|w(t)|^q e^{-t}dt.
\end{equation}

Next, we consider the case when $t\in[T_0,\infty)$. Applying (\ref{equ3.10}),  we have
\begin{align}
&\int^{\infty}_{T_0}e^{\beta|w(t)|^{N'}}|w(t)|^q e^{-t}dt \notag\\
\le&\int^{\infty}_{T_0}\exp[(\beta(1+\epsilon)-1)(t-T_0)+\beta C_{\epsilon}-T_0][(1+\epsilon)(t-T_0)+C_{\epsilon}]^{q-\frac{q}{N}}dt \notag\\
=&I_1+J_1=I_1+I_2+J_2=...=J_{N_q}+\sum^{N_q}_{k=1}I_k, \notag
\end{align}
where $N_q=\lceil q-\frac{q}{N}\rceil$. Then from the calculation we know
\begin{equation}\label{equ3.22}
I_k\le\frac{(N_q-1)!(1+\epsilon)^{k-1}C^{N_q-k}_{\epsilon}e^{\beta C_{\epsilon}}}{(N_q-k)!(1-\beta(1+\epsilon))^k}e^{-T_0},\quad\quad k=1,2,...,N_q
\end{equation}

Then we estimate $J_{N_q}$ and get
\begin{align}
J_{N_q}\le\frac{(N_q-1)!(1+\epsilon)^k}{(N_q-k-1)!(1-\beta(1+\epsilon))^k}\int^{\infty}_{T_0}A_{\beta}(t)\big((1+\epsilon)(t-T_0)+C_{\epsilon}\big)^{q-\frac{q}{N}-N_q}dt. \notag
\end{align}
Noticing  $q-\frac{q}{N}-N_q\in(-1,0]$, thus
\begin{displaymath}
\big((1+\epsilon)(t-T_0)+C_{\epsilon}\big)^{q-\frac{q}{N}-N_q}\le C^{q-\frac{q}{N}-N_q}_{\epsilon},
\end{displaymath}
therefore
\begin{equation}\label{equ3.23}
J_{N_q}\le\frac{(N_q-1)!(1+\epsilon)^k}{(N_q-k-1)!(1-\beta(1+\epsilon))^k}\frac{e^{\beta C_{\epsilon}}C^{q-\frac{q}{N}-N_q}_{\epsilon}}{1-\beta(1+\epsilon)}e^{-T_0}.
\end{equation}

Since $\epsilon$ only depends  on $\beta$, we choose $C_{2,\beta}$ as follows
\begin{displaymath}
C_{2,\beta}=\frac{(N_q-1)!(1+\epsilon)^k}{(N_q-k-1)!(1-\beta(1+\epsilon))^k}\frac{e^{\beta C_{\epsilon}}C^{q-\frac{q}{N}-N_q}_{\epsilon}}{1-\beta(1+\epsilon)}
+\sum^{N}_{k=1}\frac{(N_q-1)!(1+\epsilon)^{k-1}C^{N_q-k}_{\epsilon}e^{\beta C_{\epsilon}}}{(N_q-k)!(1-\beta(1+\epsilon))^k}.
\end{displaymath}
Thus, we can conclude
\begin{equation}\label{equ3.24}
\int^{\infty}_{T_0}e^{\beta|w(t)|^{N'}}|w(t)|^N e^{-t}dt\le C_{2,\beta}\int^{\infty}_{T_0}|w(t)|^N e^{-t}dt.
\end{equation}

Combining (\ref{equ3.21}) with (\ref{equ3.24}) together,  we have then proved Lemma 3.2.
\end{proof}

\section{Proofs of Theorem 1.2 and Theorem 1.3.}

Now we begin to consider Theorem 1.2. As in Section 2, it suffices for us to prove the inequality of the special case $s=t$, which states that, under the assumption (\ref{equ1.1}), there exists a positive constant $C=C(N,t,\alpha)$ such that the inequality
\begin{equation}\label{equ4.1}
\int_{\mathbb{R}^N}e^{\alpha|u|^{N'}}|u|^N\frac{dx}{|x|^t}\le C\int_{\mathbb{R}^N}\frac{|u|^N}{|x|^t}dx
\end{equation}
holds for all functions $u\in X^{1,N}_s$ with $\|\nabla u\|_{L^N(\mathbb{R}^N)}\le 1$.\\

 Once we have proved the special case (\ref{equ4.1}), the general case $s<t$ follows immediately by applying the Caffarelli-Kohn-Nirenberg inequality (\ref{equ2.6}).\\

However, since the functions $u$ under consideration are not required to be spherically symmetric, we cannot use the symmetrization method
to reduce the proof of (\ref{equ4.1}) to only spherically symmetric functions due to the existence of  the weight $\frac{1}{|x|^t}$.
Therefore, the method  used in \cite{INW} does not work here. To overcome this difficulty, we will develop a new argument of change of variables to attack this problem. \\

Now we begin the proof of Theorem 1.2.
\begin{proof}
Let $0<\alpha<\alpha_{N,t}$ and let $u\in X^{1,N}_s$ with $\|\nabla u\|_{L^N(\mathbb{R}^N)}\le 1$. We define the function $v\in W^{1,N}(\mathbb{R}^N)$ for $x\in\mathbb{R}^N$ in the same way  as in (\ref{equ2.7}),
\begin{displaymath}
v(x):=\Big(\frac{N-t}{N}\Big)^{\frac{1}{N'}}u(|x|^{\frac{t}{N-t}}x).
\end{displaymath}

Form (\ref{equ2.13}),(\ref{equ2.14}), we have
\begin{equation}\label{equ4.2}
\|u\|_{L^N(\mathbb{R}^N;|x|^{-t}dx)}=\frac{N}{N-t}\|v\|_{L^N(\mathbb{R}^N)},
\end{equation}
\begin{equation}\label{equ4.3}
\|\nabla u\|_{L^N(\mathbb{R}^N)}\ge\|\nabla v\|_{L^N(\mathbb{R}^N)}.
\end{equation}

Direct computations show
\begin{align}\label{equ4.4}
&\int_{\mathbb{R}^N}e^{\alpha|u(y)|^{\frac{N}{N-1}}}|u(y)|^N\frac{dy}{|y|^t} \notag\\
=&\sum^{\infty}_{i=0}\int_{\mathbb{R}^N}\frac{(\alpha|u(y)|^{N'})^i}{i!}|u(y)|^{N'(N-1)}\frac{dy}{|y|^t} \notag\\
=&\sum^{\infty}_{i=0}\int_{\mathbb{R}^N}\frac{N}{N-t}\frac{\alpha^i}{i!}\big|u(|x|^{\frac{t}{N-t}}x)\big|^{N'(i+N-1)}dx \notag\\
=&\sum^{\infty}_{i=0}\int_{\mathbb{R}^N}\Big(\frac{N}{N-t}\Big)^{i+N}\frac{\alpha^i}{i!}|v(x)|^{N'(i+N-1)}dx \notag\\
=&\Big(\frac{N}{N-t}\Big)^N\sum^{\infty}_{i=0}\int_{\mathbb{R}^N}\frac{(\frac{N}{N-t}\alpha|v(x)|^{N'})^i}{i!}|v(x)|^N dx \notag\\
=&\Big(\frac{N}{N-t}\Big)^N\int_{\mathbb{R}^N}\exp\big(\frac{N}{N-t}\alpha|v(x)|^{N'}\big)|v(x)|^N dx.
\end{align}

Since $0<\frac{N}{N-t}\alpha<\frac{N}{N-t}\alpha_{N,t}=N\omega^{\frac{1}{N-1}}_{N-1}$ and $\|\nabla v\|_{L^N(\mathbb{R}^N)}\le 1$, by applying Lemma 3.1 and (\ref{equ4.2}), (\ref{equ4.4}) we have
\begin{align}
\int_{\mathbb{R}^N}e^{\alpha|u|^{N'}}|u|^N\frac{dx}{|x|^t} \notag=&(\frac{N}{N-t})^N\int_{\mathbb{R}^N}e^{\frac{N}{N-t}\alpha|v|^{N'}}|v|^N dx \notag\\
\le&(\frac{N}{N-t})^N C_{\alpha}\|v\|^N_{L^N(\mathbb{R}^N)} \notag\\
=&C_{\alpha}\|u\|^N_{L^N(\mathbb{R}^N;|x|^{-t}dx)}.
\end{align}
Thus,  we have proved  inequality (\ref{equ1.3}) in Theorem 1.2.\\

Next we want to show $\alpha_{N,t}$ is the sharp constant.\\

Applying the same test sequence (\ref{equ2.17}) again, we have $\|\nabla u_k\|_{L^N(\mathbb{R}^N)}=1$ for all $k\in\mathbb{N}$ and $\|u_k\|_{L^N(\mathbb{R}^N;|x|^{-s}dx)}=o(1)$ as $k\to\infty$. By direct calculations, we have
\begin{align}
&\int_{\mathbb{R}^N}e^{\alpha_{N,t}|u_k|^{N'}}|u_k|^N\frac{dx}{|x|^t} \notag\\
\ge&\int^{e^{-\frac{k}{N-t}}}_{0}\omega_{N-1}e^k\frac{1}{\omega_{N-1}}(\frac{k}{N-t})^{N-1}r^{N-t-1}dr \notag\\
=&\frac{k^{N-1}}{(N-t)^N}, \notag
\end{align}
thus we have
\begin{align}
\frac{\int_{\mathbb{R}^N}e^{\alpha_{N,t}|u_k|^{N'}}|u_k|^N\frac{dx}{|x|^t}}{\|u_k\|^{\frac{N(N-t)}{N-s}}_{L^N(\mathbb{R}^N;|x|^{-s}dx)}}\to\infty \quad\quad as\,\,k\to\infty, \notag
\end{align}
which implies inequality (\ref{equ1.3}) fails when $\alpha=\alpha_{N,t}$, hence we have finished the proof of Theorem 1.2.
\end{proof}

We now start to prove Theorem 1.3.
The method of proving Theorem 1.3 is similar to Theorem 1.2. By using the Caffarelli-Kohn-Nirenberg inequality (\ref{equ2.6}) it suffices for us to prove the special case $s=t$, which states  the following inequality
\begin{equation}\label{equ4.6}
\int_{\mathbb{R}^N}e^{\alpha|u|^{N/(N-1)}}|u|^q\frac{dx}{|x|^t}\le\int_{\mathbb{R}^N}\frac{|u|^q}{|x|^t}dx
\end{equation}
holds for all functions $u\in L^q(\mathbb{R}^N;|x|^{-t}ds)\cap \dot{W}^{1,N}(\mathbb{R})$ with $\|\nabla u\|_{L^N(\mathbb{R}^N)}\le 1$.

\begin{proof}
Let $0<\alpha<\alpha_{N,t}$ and let $u\in L^q(\mathbb{R}^N;|x|^{-t}dx)\cap \dot{W}^{1,q}(\mathbb{R}^N)$ with $\|\nabla u\|_{L^N(\mathbb{R}^N)}\le 1$. We define the function $v\in W^{1,q}(\mathbb{R}^N)$   in the same way as in (\ref{equ2.7}),
\begin{displaymath}
v(x):=\Big(\frac{N-t}{N}\Big)^{\frac{N-1}{N}}u(|x|^{\frac{t}{N-t}}x).
\end{displaymath}

For $j\ge q$, applying the change of variables $y=|x|^{\frac{t}{N-t}}x$, we have,
\begin{align}\label{equ4.7}
\int_{\mathbb{R}^N}|v(x)|^j dx=&\int_{\mathbb{R}^N}\Big(\frac{N-t}{N}\Big)^{\frac{j}{N'}}\big|u(|x|^{\frac{t}{N-t}}x)\big|^j dx \notag\\
=&\Big(\frac{N-t}{N}\Big)^{\frac{j}{N'}+1}\int_{\mathbb{R}^N}\frac{|u(y)|^j}{|y|^t}dy,
\end{align}
similarly we can  have
\begin{align}\label{equ4.8}
\int_{\mathbb{R}^N}e^{\alpha|u(y)|^{N'}}|u(y)|^q\frac{dy}{|y|^t}=&\sum^{\infty}_{i=0}\int_{\mathbb{R}^N}\frac{(\alpha|u(y)|^{N'})^i}{i!}|u(y)|^q\frac{dy}{|y|^t} \notag\\
=&\sum^{\infty}_{i=0}\int_{\mathbb{R}^N}\frac{\alpha^i}{i!}\Big(\frac{N}{N-t}\Big)^{i+\frac{q}{N'}+1}
\Big(\frac{N-t}{N}\Big)^{i+\frac{q}{N'}+1}\big|u(y)\big|^{N'i+q}\frac{dy}{|y|^t} \notag\\
=&\Big(\frac{N}{N-t}\Big)^{\frac{q}{N'}+1}\sum^{\infty}_{i=0}\int_{\mathbb{R}^N}\frac{\alpha^i}{i!}\Big(\frac{N}{N-t}\Big)^i|v(x)|^{N'i+q}dx \notag\\
=&\Big(\frac{N}{N-t}\Big)^{\frac{q}{N'}+1}\int_{\mathbb{R}^N}\exp\big(\frac{N}{N-t}\alpha|v(x)|^{N'}\big)|v(x)|^q dx.
\end{align}

By (\ref{equ4.3}),  we have
\begin{align}
\|\nabla v\|^N_{L^N(\mathbb{R}^N)}\le\|\nabla u\|^N_{L^N(\mathbb{R}^N)}\le 1.
\end{align}

Since $0<\frac{N}{N-t}\alpha<\frac{N}{N-t}\alpha_{N,t}=N\omega^{\frac{1}{N-1}}_{N-1}$, and $\|\nabla v\|^N_{L^N(\mathbb{R}^N)}\le 1$, we apply (\ref{equ4.7}), (\ref{equ4.8}) and Lemma 3.2 to get
\begin{align}\label{equ4.10}
\int_{\mathbb{R}^N}e^{\alpha|u|^{N'}}|u|^q\frac{dx}{|x|^t} \notag=&\Big(\frac{N}{N-t}\Big)^{\frac{q}{N'}+1}\int_{\mathbb{R}^N}e^{\frac{N}{N-t}\alpha|v|^{N'}}|v|^q dx \notag\\
\le&\Big(\frac{N}{N-t}\Big)^{\frac{q}{N'}+1} C\|v\|^q_{L^q(\mathbb{R}^N)} \notag\\
=&C\|u\|^q_{L^q(\mathbb{R}^N;|x|^{-t}dx)},
\end{align}
where $C=C(N,\alpha,q,t)$.\\

Therefore (\ref{equ4.6}) has been established and we have proved inequality (\ref{equ1.4}).\\

Next,  we will show $\alpha_{N,t}$ is the sharp constant for our inequality.\\

For $k\in\mathbb{N}$, define a sequence $u_k$ of radially symmetric function in $u\in L^N(\mathbb{R}^N;|x|^{-s}ds)\cap \dot{W}^{1,N}(\mathbb{R})$ by
\begin{align}\label{equ4.11}
u_k(x)=\left\{
\begin{aligned}
&0 & \text{if } |x|\ge 1 \\
&\Big(\frac{N-t}{\omega_{N-1}k}\Big)^{\frac{1}{N}}\Big(\frac{q-t}{N-t}\Big)\log{\frac{1}{|x|}}& \text{if } e^{-\frac{k}{q-t}}<|x|<1 \\
&\Big(\frac{1}{\omega_{N-1}}\Big)^{\frac{1}{N}}\Big(\frac{k}{N-t}\Big)^{\frac{1}{N'}}& \text{if } 0\le|x|\le e^{-\frac{k}{q-t}}
\end{aligned}
\right..
\end{align}

Direct computations show that $\|\nabla u_k\|_{L^N(\mathbb{R}^N)}=1$ for all $k\in\mathbb{N}$ and $\|u_k\|_{L^q(\mathbb{R}^N;|x|^{-s}dx)}=o(1)$ as $k\to\infty$.

Moreover,  we have
\begin{align}
&\int_{\mathbb{R}^N}e^{\alpha_{N,t}|u_k|^{N'}}|u_k|^q\frac{dx}{|x|^t} \notag\\
\ge&\int^{e^{-\frac{k}{q-t}}}_{0}\omega_{N-1}e^k\Big(\frac{1}{\omega_{N-1}}\Big)^{\frac{q}{N}}\Big(\frac{k}{N-t}\Big)^{\frac{q}{N'}}r^{q-t-1}dr \notag\\
=&\frac{\omega^{1-\frac{q}{N}}_{N-1}}{q-t}\Big(\frac{k}{N-t}\Big)^{\frac{q}{N'}}, \notag
\end{align}
thus we have
\begin{align}
\frac{\int_{\mathbb{R}^N}e^{\alpha_{N,t}|u_k|^{N'}}|u_k|^q\frac{dx}{|x|^t}}{\|u_k\|^{\frac{N(N-t)}{N-s}}_{L^q(\mathbb{R}^N;|x|^{-s}dx)}}\to\infty \quad\quad as\,\,k\to\infty, \notag
\end{align}
which shows that $\alpha_{N,t}$ is the sharp constant.
\end{proof}

\section{Proof of Theorem 1.4.}

We first introduce the rearrangement function for a measurable function $u$ on $\mathbb{R}^N$.  We define the distribution function of $u$ by
\begin{displaymath}
a_u(\lambda):=|\{x\in\mathbb{R}^N;|u(x)|>\lambda\}|,
\end{displaymath}
then the rearrangement function $u^\sharp:\mathbb{R}^N\to[0,\infty]$ is defined by
\begin{align}
\left\{
\begin{array}{rcl}
u^\sharp(0)&=&\text{ess.}\sup(u), \\
u^\sharp(r)&=&\inf\{\lambda|a_u(\lambda)<r\},\quad r>0. \notag
\end{array}
\right.
\end{align}

To prove the first part of Theorem 1.4,  let us recall some Lemmas and one Corollary proved in \cite{INW}. The first one is well-known and follows easily from the Hardy-Littlewood inequality. \\
\\
\textbf{Lemma 5.1}  \emph{Let $N\ge2$, $0\le t<N$ and $q\ge N$. Then it holds
\begin{equation}\label{equ5.1}
\|u\|_{L^q(\mathbb{R}^N;|x|^{-t}dx)}\le\|u^\sharp\|_{L^q(\mathbb{R}^N;|x|^{-t}dx)}
\end{equation}
for all functions $u$ so that $u^\sharp\in L^q(\mathbb{R}^N;|x|^{-t}dx)$.}\\
\\
\textbf{Lemma 5.2} (Lemma A.3 in \cite{INW}) \emph{Let $N\ge2$ and let (s,t,q) be exponents satisfying either
\begin{displaymath}
-\infty<s<t<N \text{ and } N\le q<\infty\quad \text{ or }\quad -\infty<s=t<N \text{ and } N<q<\infty
\end{displaymath}
Then the embedding
\begin{displaymath}
X^{1,N}_{s,rad}\hookrightarrow L^q(\mathbb{R}^N;|x|^{-t}dx)
\end{displaymath}
is compact.}\\
\\
\textbf{Proposition 5.3} (Corollary 1.4 in \cite{INW}) \emph{Assume (\ref{equ1.1}) with $s=0$, then
\begin{displaymath}
\mu_{N,0,t,\alpha}(\mathbb{R}^N)=\sup_{\substack{u\in H^{1,N}(\mathbb{R}^N)\\\|\nabla u\|_{L^N(\mathbb{R}^N)}=1}}
\frac{\int_{\mathbb{R}^N}\Phi_N(\alpha|u|^{N'})\frac{dx}{|x|^t}}{\|u\|^{N-t}_{L^N(\mathbb{R}^N)}}
\end{displaymath}
is attained.}\\

Now we are in the position to prove the first part of Theorem 1.4.

\begin{proof}[Proof of Theorem 1.4 (i)]
Consider in (\ref{equ1.5}) we have
\begin{equation}\label{equ5.2}
F_{N,s,t,\alpha}(u)=\frac{\int_{\mathbb{R}^N}\Phi_N(\alpha|u|^{N'})\frac{dx}{|x|^t}}{\|u\|^{\frac{N(N-t)}{N-s}}_{L^N(\mathbb{R}^N;|x|^{-s}dx)}},
\end{equation}
define a new function $v\in H^{1,N}(\mathbb{R}^N)$ as following,
\begin{equation}\label{equ5.3}
v(x):=\Big(\frac{N-s}{N}\Big)^{\frac{1}{N'}}u(|x|^{\frac{s}{N-s}}x).
\end{equation}
similar to (\ref{equ2.13}) we can get,
\begin{equation}\label{equ5.4}
\|u\|_{L^N(\mathbb{R}^N;|x|^{-s}dx)}=\frac{N}{N-s}\|v\|_{L^N(\mathbb{R}^N)},
\end{equation}
the direct computation also show us
\begin{equation}\label{equ5.5}
\int_{\mathbb{R}^N}\Phi_N(\alpha|u|^{N'})\frac{dx}{|x|^t}=\frac{N}{N-s}\int_{\mathbb{R}^N}\Phi_N(\alpha\frac{N}{N-s}|v|^{N'})\frac{dx}{|x|^{\frac{N(t-s)}{N-s}}}.
\end{equation}

Let $t^*=\frac{N(t-s)}{N-s}\in[0,N)$, we notice $0<\alpha\frac{N}{N-s}<\alpha_{N,t^*}$. Therefore according to (\ref{equ5.2}), (\ref{equ5.4}), (\ref{equ5.5}) and Proposition 5.3 we have
\begin{align}
\mu_{N,s,t,\alpha}(\mathbb{R}^N)&=\sup_{\substack{u\in X^{1,N}_s\\\|\nabla u\|_{L^N(\mathbb{R}^N)}=1}}F_{N,s,t,\alpha}(u) \notag\\
&=\Big(\frac{N}{N-s}\Big)^{1+t^*-N}\Bigg(\sup_{\substack{v\in H^{1,N}(\mathbb{R}^N)\\\|\nabla v\|_{L^N(\mathbb{R}^N)}=1}}
\frac{\int_{\mathbb{R}^N}\Phi_N(\alpha\frac{N}{N-s}|v|^{N'})\frac{dx}{|x|^{t^*}}}{\|v\|^{N-t^*}_{L^N(\mathbb{R}^N)}}\Bigg). \notag
\end{align}
is attained. Hence we proved part (i) of Theorem 1.4.
\end{proof}

In order to prove the second part of Theorem 1.4, we want to show $\nu_{N,s,t,\alpha}$ is attained when the functions $u$ are radially symmetric.\\
\\
\textbf{Proposition 5.4} \emph{Assume (\ref{equ1.1}) holds, and let $G_{N,s,t,\alpha}(u)$ be defined as in (\ref{equ1.6}).  Then
\begin{displaymath}
\nu_{N,s,t,\alpha,rad}(\mathbb{R}^N)=\sup_{\substack{u\in X^{1,N}_{s,rad}\\\|\nabla u\|_{L^N(\mathbb{R}^N)}=1}}G_{N,s,t,\alpha}(u)
\end{displaymath}
is attained.}\\

From Lemma 5.2 we notice the non-compactness for embedding $X^{1,N}_{s,rad}\hookrightarrow L^N(\mathbb{R}^N;|x|^{-t}dx)$ when $s=t$, hence we establish the following lemma first.\\
\\
\textbf{Lemma 5.5} \emph{Assume (\ref{equ1.1}), and let $\{u_n\}$ be a bounded sequence that belongs to $X^{1,N}_{s,rad}$ with $\|\nabla u\|_{L^N(\mathbb{R}^N)}=1$. Also we have
\begin{displaymath}
u_n\rightharpoonup u\quad\text{ weakly in}\quad X^{1,N}_{s,rad}
\end{displaymath}
as $n\to\infty$, then the following convergence holds as $n\to\infty$.
\begin{equation}\label{equ5.6}
\int_{\mathbb{R}^N}(e^{\alpha|u_n|^{N'}}|u_n|^N-|u_n|^N)\frac{dx}{|x|^t}\to\int_{\mathbb{R}^N}(e^{\alpha|u|^{N'}}|u|^N-|u|^N)\frac{dx}{|x|^t}.
\end{equation}}

We remark here that a similar lemma when we replace $e^{\alpha|u_n|^{N'}}|u_n|^N-|u_n|^N$ by $\Phi_N(\alpha|u_n|^{N'})-\frac{\alpha^{N-1}}{(N-1)!}|u_n|^{N-1}$ was carried out in \cite{INW}, and such an idea appears in a number of works, e.g., \cite{CC}  and \cite{LR}, etc. We include a proof for our case here for the sake of completeness.

\begin{proof}
Let $\Psi_k(\tau):=\sum^{\infty}_{i=k}\frac{\alpha^{(i-(N-1))}}{(i-(N-1))!}|\tau|^{N'i}$, where $N\ge 2$, $N-1\le k\in\mathbb{N}$ and $0<\alpha<\alpha_{N,t}$.

It is not hard to check $e^{\alpha|u|^{N'}}|u|^N=\Psi_{N-1}(|u|)$ and $e^{\alpha|u|^{N'}}|u|^N-|u|^N=\Psi_{N}(|u|)$, therefore (\ref{equ5.6}) becomes
\begin{displaymath}
\int_{\mathbb{R}^N}\Psi_N(|u_n|)\frac{dx}{|x|^t}\to\int_{\mathbb{R}^N}\Psi_N(|u|)\frac{dx}{|x|^t}
\end{displaymath}
as $n\to\infty$.

Direct calculation show us,
\begin{align}\label{equ5.7}
\Psi'_N(\tau)&=N'\alpha\tau^{\frac{1}{N-1}}\sum^{\infty}_{i=N-1}\frac{\alpha^{(i-(N-1))}(i+1)}{(i-(N-2))!}|\tau|^{N'i} \notag\\
&\le NN'\alpha\tau^{\frac{1}{N-1}}\sum^{\infty}_{i=N-1}\frac{\alpha^{(i-(N-1))}}{(i-(N-1))!}|\tau|^{N'i} \notag\\
&=NN'\alpha\tau^{\frac{1}{N-1}}\Psi_{N-1}(\tau),
\end{align}
by the mean value theorem and the convexity of $\Psi_{N-1}$ we know there exists some $\theta\in[0,1]$ such that
\begin{align}\label{equ5.8}
&|\Psi_N(|u_n|)-\Psi_N(|u|)| \notag\\
\le&\Psi'_N(\theta|u_n|+(1-\theta)|u|)|u_n-u| \notag\\
\le&NN'\alpha(\theta|u_n|+(1-\theta)|u|)^{\frac{1}{N-1}}\Psi_{N-1}(\theta|u_n|+(1-\theta)|u|)|u_n-u| \notag\\
\le&NN'\alpha(\theta|u_n|+(1-\theta)|u|)^{\frac{1}{N-1}}(\theta\Psi_{N-1}(|u_n|)+(1-\theta)\Psi_{N-1}(|u|))|u_n-u| \notag\\
\le&NN'\alpha(|u_n|+|u|)^{\frac{1}{N-1}}(\Psi_{N-1}(|u_n|)+\Psi_{N-1}(|u|))|u_n-u|.
\end{align}

Then take the numbers $a,b,c>1$ satisfy $\frac{1}{a}+\frac{1}{b}+\frac{1}{c}=1$ which we will choose later, by Holder inequality we have
\begin{align}\label{equ5.9}
&\Big{|}\int_{\mathbb{R}^N}(\Psi_N(|u_n|)-\Psi_N(|u|))\frac{dx}{|x|^t}\Big{|} \notag\\
\le&NN'\alpha\int_{\mathbb{R}^N}(|u_n|+|u|)^{\frac{1}{N-1}}(\Psi_{N-1}(|u_n|)+\Psi_{N-1}(|u|))|u_n-u|\frac{dx}{|x|^t} \notag\\
\le&NN'\alpha\big{\|}|u_n|+|u|\big{\|}^{\frac{1}{N-1}}_{L^{\frac{a}{N-1}}(\mathbb{R}^N;|x|^{-t}dx)}
\big{\|}\Psi_{N-1}(|u_n|)+\Psi_{N-1}(|u|)\big{\|}_{L^b(\mathbb{R}^N;|x|^{-t}dx)}\big{\|}u_n-u\big{\|}_{L^c(\mathbb{R}^N;|x|^{-t}dx)}.
\end{align}

From Caffarelli-Kohn-Nirenberg inequality (\ref{equ2.6}) we obtain the boundedness of $X^{1,N}_s\hookrightarrow L^\frac{a}{N-1}(\mathbb{R}^N;|x|^{-t}dx)$ for $\frac{a}{N-1}\ge N$, which gives us
\begin{equation}\label{equ5.10}
\big{\|}|u_n|+|u|\big{\|}_{L^{\frac{a}{N-1}}(\mathbb{R}^N;|x|^{-t}dx)}\le C(\|u_n\|_{X^{1,n}_s}+\|u\|_{X^{1,n}_s})\le C,
\end{equation}

And we could choose $b>1$ sufficiently close to 1 such that $b\alpha<\alpha_{N,t}$, from Lemma 5.2 and $bN>N$ we know $\|u_n\|_{L^{bN}(\mathbb{R}^N;|x|^{-s}dx)}\le\|u_n\|_{L^N(\mathbb{R}^N;|x|^{-s}dx)}$, combine with Theorem 1.3 we have
\begin{align}\label{equ5.11}
\big{\|}\Psi_{N-1}(|u_n|)\big{\|}_{L^b(\mathbb{R}^N;|x|^{-t}dx)}=&\big{(}\int_{\mathbb{R}^N}(e^{\alpha|u_n|^{N'}}|u_n|^N)^b\frac{dx}{|x|^t}\big{)}^{\frac{1}{b}} \notag\\
=&\big{(}\int_{\mathbb{R}^N}e^{b\alpha|u_n|^{N'}}|u_n|^{bN}\frac{dx}{|x|^t}\big{)}^{\frac{1}{b}} \notag\\
\le&C\|u_n\|^{\frac{N(N-t)}{b(N-s)}}_{L^{bN}(\mathbb{R}^N;|x|^{-s}dx)} \notag\\
\le&C\|u_n\|^{\frac{N(N-t)}{b(N-s)}}_{L^N(\mathbb{R}^N;|x|^{-s}dx)}\le C,
\end{align}

similarly we obtain
\begin{equation}\label{equ5.12}
\big{\|}\Psi_{N-1}(|u|)\big{\|}_{L^b(\mathbb{R}^N;|x|^{-t}dx)}\le C,
\end{equation}

Furthermore, from Lemma 5.2 we have the compactness for $X^{1,N}_{s,rad}\hookrightarrow L^c(\mathbb{R}^N;|x|^{-t}dx)$ for $c>N$, hence we have the convergence
\begin{equation}\label{equ5.13}
\big{\|}u_n-u\big{\|}_{L^c(\mathbb{R}^N;|x|^{-t}dx)}\to 0 \quad\quad as\,\, n\to\infty.
\end{equation}
combine (\ref{equ5.10})-(\ref{equ5.13}) we have
\begin{equation}\label{equ5.14}
\Big{|}\int_{\mathbb{R}^N}(\Psi_N(|u_n|)-\Psi_N(|u|))\frac{dx}{|x|^t}\Big{|}\to 0\quad\quad as\,\, n\to\infty.
\end{equation}

Therefore we proved this Lemma.
\end{proof}

Now we are in position to prove Proposition 5.4 by applying Lemma 5.5.
\begin{proof}[Proof of Proposition 5.4]
Let $\{u_n\}$ be a maximizing sequence for $\nu_{N,s,t,\alpha.rad}(\mathbb{R}^N)$, which gives us $G_{N,s,t,\alpha}(u_n)\to\nu_{N,s,t,\alpha,rad}(\mathbb{R}^N)$ as $n\to\infty$. We define a new sequence $\{v_n\}$ by $v_n(x):=u_n(\|u_n\|^{\frac{N}{N-s}}_{L^N(\mathbb{R}^N;|x|^{-s}dx)}x)$ for $x\in\mathbb{R}^N$. Then direct calculation show us
\begin{displaymath}
\|\nabla v_n\|_{L^N(\mathbb{R}^N)}=\|\nabla u_n\|_{L^N(\mathbb{R}^N)}=1, \quad \|v_n\|_{L^N(\mathbb{R}^N;|x|^{-s}dx)}=1,
\end{displaymath}
and
\begin{displaymath}
G_{N,s,t,\alpha}(v_n)=G_{N,s,t,\alpha}(u_n)\to\nu_{N,s,t,\alpha,rad}(\mathbb{R}^N) \quad\text{as }n\to\infty.
\end{displaymath}

Thus $\{v_n\}$ is also a maximizing sequence for $\nu_{N,s,t,\alpha,rad}(\mathbb{R}^N)$. Therefore, up to a subsequence, $v_n$ converges to some $v$ weakly in $X^{1,N}_{s,rad}$, then $v$ satisfies
\begin{equation}\label{equ5.15}
\max{\{\|v\|_{L^N(\mathbb{R}^N;|x|^{-s}dx)},\|\nabla v\|_{L^N(\mathbb{R}^N)}\}}\le1.
\end{equation}

First we consider the case when $s=t$, we could assume $\nu_{N,s,\alpha,rad}(\mathbb{R}^N)=\nu_{N,s,s,\alpha,rad}(\mathbb{R}^N)$ and $G_{N,s,\alpha}(u)=G_{N,s,s,\alpha}(u)$ for $u\in X^{1,N}_{s,rad}$.\\

Apply Lemma 5.5 and let $n\to\infty$ we see that,
\begin{align}\label{equ5.16}
\nu_{N,s,\alpha,rad}(\mathbb{R}^N)&=G_{N,s,\alpha}(v_n)+o(1) \notag\\
&=\int_{\mathbb{R}^N}e^{\alpha|v_n|^{N'}}|v_n|^N\frac{dx}{|x|^s}+o(1) \notag\\
&=1+\int_{\mathbb{R}^N}(e^{\alpha|v_n|^{N'}}|v_n|^N-|v_n|^N)\frac{dx}{|x|^s}+o(1) \notag\\
&=1+\int_{\mathbb{R}^N}(e^{\alpha|v|^{N'}}|v|^N-|v|^N)\frac{dx}{|x|^s}.
\end{align}

Pick up any $u_0\in X^{1,N}_{s,rad}$ satisfying $\|\nabla u_0\|_{L^N(\mathbb{R}^N)}=1$ we could see,
\begin{align}
\nu_{N,s,\alpha,rad}(\mathbb{R}^N)&\ge G_{N,s,\alpha}(u_0) \notag\\
&=\frac{\int_{\mathbb{R}^N}e^{\alpha|u_0|^{N'}}|u_0|^N\frac{dx}{|x|^s}}{\|u_0\|^N_{L^N(\mathbb{R}^N;|x|^{-s}dx)}} \notag\\
&=\frac{\sum^{\infty}_{j=N-1}\frac{\alpha^{j-(N-1)}}{(j-(N-1))!}\|u_0\|^{N'j}_{L^{N'j}(\mathbb{R}^N;|x|^{-s}dx)}}{\|u_0\|^N_{L^N(\mathbb{R}^N;|x|^{-s}dx)}} \notag\\
&=1+\frac{\sum^{\infty}_{j=N}\frac{\alpha^{j-(N-1)}}{(j-(N-1))!}\|u_0\|^{N'j}_{L^{N'j}(\mathbb{R}^N;|x|^{-s}dx)}}{\|u_0\|^N_{L^N(\mathbb{R}^N;|x|^{-s}dx)}}>1, \notag
\end{align}
combine with (\ref{equ5.16}) we know the $\int_{\mathbb{R}^N}(e^{\alpha|v|^{N'}}|v|^N-|v|^N)\frac{dx}{|x|^s}>0$, which implies $v$ is not identity 0. Since $\|v\|_{L^N(\mathbb{R}^N;|x|^{-s}dx)}\le1$, we get
\begin{align}\label{equ5.17}
\nu_{N,s,\alpha,rad}(\mathbb{R}^N)&\le 1+\frac{\int_{\mathbb{R}^N}(e^{\alpha|v|^{N'}}|v|^N-|v|^N)\frac{dx}{|x|^s}}{\|v\|^N_{L^N(\mathbb{R}^N;|x|^{-s}dx)}} \notag\\
&=\frac{\int_{\mathbb{R}^N}(e^{\alpha|v|^{N'}}|v|^N)\frac{dx}{|x|^s}}{\|v\|^N_{L^N(\mathbb{R}^N;|x|^{-s}dx)}}=G_{N,s,\alpha}(v).
\end{align}

Therefore if we can prove $\|\nabla v\|_{L^N(\mathbb{R}^N)}=1$ the theorem proved. Since we already know $\|\nabla v\|_{L^N(\mathbb{R}^N)}\le1$, it suffices to show $\|\nabla v\|_{L^N(\mathbb{R}^N)}\ge1$. So we have
\begin{align}\label{equ5.18}
\nu_{N,s,\alpha,rad}(\mathbb{R}^N)&\ge G_{N,s,\alpha}(\frac{v}{\|\nabla v\|_{L^N(\mathbb{R}^N)}}) \notag\\
&=\frac{\|\nabla v\|^N_{L^N(\mathbb{R}^N)}}{\|v\|^N_{L^N(\mathbb{R}^N;|x|^{-s}dx)}}
\int_{\mathbb{R}^N}e^{\alpha|\frac{v}{\|\nabla v\|_{L^N(\mathbb{R}^N)}}|^{N'}}|\frac{v}{\|\nabla v\|_{L^N(\mathbb{R}^N)}}|^N\frac{dx}{|x|^s} \notag\\
&=\frac{\|\nabla v\|^N_{L^N(\mathbb{R}^N)}}{\|v\|^N_{L^N(\mathbb{R}^N;|x|^{-s}dx)}}
\sum^{\infty}_{j=N-1}\frac{\alpha^{j-(N-1)}}{(j-(N-1))!}\frac{\|v\|^{N'j}_{L^{N'j}(\mathbb{R}^N;|x|^{-s}dx)}}{\|\nabla v\|^{N'j}_{L^N(\mathbb{R}^N)}} \notag\\
&=\sum^{\infty}_{j=N-1}\frac{\alpha^{j-(N-1)}}{(j-(N-1))!}\frac{\|v\|^{N'j}_{L^{N'j}(\mathbb{R}^N;|x|^{-s}dx)}}{\|v\|^N_{L^N(\mathbb{R}^N;|x|^{-s}dx)}}
\|\nabla v\|^{N-N'j}_{L^N(\mathbb{R}^N)} \notag\\
&\ge 1+\alpha\frac{\|v\|^{N'N}_{L^{N'N}(\mathbb{R}^N;|x|^{-s}dx)}}{\|v\|^N_{L^N(\mathbb{R}^N;|x|^{-s}dx)}}\|\nabla v\|^{-N'}_{L^N(\mathbb{R}^N)}
+\sum^{\infty}_{j=N+1}\frac{\alpha^{j-(N-1)}}{(j-(N-1))!}\frac{\|v\|^{N'j}_{L^{N'j}(\mathbb{R}^N;|x|^{-s}dx)}}{\|v\|^N_{L^N(\mathbb{R}^N;|x|^{-s}dx)}} \notag\\
&=1+\sum^{\infty}_{j=N}\frac{\alpha^{j-(N-1)}}{(j-(N-1))!}\frac{\|v\|^{N'j}_{L^{N'j}(\mathbb{R}^N;|x|^{-s}dx)}}{\|v\|^N_{L^N(\mathbb{R}^N;|x|^{-s}dx)}}
+\alpha(\frac{1}{\|\nabla v\|^{N'}_{L^N(\mathbb{R}^N)}}-1)\frac{\|v\|^{N'N}_{L^{N'N}(\mathbb{R}^N;|x|^{-s}dx)}}{\|v\|^N_{L^N(\mathbb{R}^N;|x|^{-s}dx)}} \notag\\
&=G_{N,s,\alpha}(v)+\alpha(\frac{1}{\|\nabla v\|^{N'}_{L^N(\mathbb{R}^N)}}-1)\frac{\|v\|^{N'N}_{L^{N'N}(\mathbb{R}^N;|x|^{-s}dx)}}{\|v\|^N_{L^N(\mathbb{R}^N;|x|^{-s}dx)}},
\end{align}
combine with (\ref{equ5.17}) we have
\begin{displaymath}
\nu_{N,s,\alpha,rad}(\mathbb{R}^N)\ge\nu_{N,s,\alpha,rad}(\mathbb{R}^N)
+\alpha(\frac{1}{\|\nabla v\|^{N'}_{L^N(\mathbb{R}^N)}}-1)\frac{\|v\|^{N'N}_{L^{N'N}(\mathbb{R}^N;|x|^{-s}dx)}}{\|v\|^N_{L^N(\mathbb{R}^N;|x|^{-s}dx)}},
\end{displaymath}
which directly tell us $\|\nabla v\|_{L^N(\mathbb{R}^N)}\ge1$, then it follows $\|\nabla v\|_{L^N(\mathbb{R}^N)}=1$. Hence we shows that $v$ is a maximizer for $\nu_{N,s,\alpha,rad}(\mathbb{R}^N)$.\\

Then consider the case $s<t$, by Lemma 5.2 we have the compactness of the embedding $X^{1,N}_{s,rad}\hookrightarrow L^N(\mathbb{R}^N;|x|^{-t}dx)$. Hence we have the convergence as $n\to\infty$,
\begin{align}
\nu_{N,s,\alpha,rad}(\mathbb{R}^N)&=G_{N,s,\alpha}(v_n)+o(1) \notag\\
&=\int_{\mathbb{R}^N}e^{\alpha|v_n|^{N'}}|v_n|^N\frac{dx}{|x|^s}+o(1) \notag\\
&=\int_{\mathbb{R}^N}e^{\alpha|v|^{N'}}|v|^N\frac{dx}{|x|^s}, \notag
\end{align}
which implies $v$ is not identity 0. Then the following we could using the same method as we used when $s=t$ to prove $\|\nabla v\|_{L^N(\mathbb{R}^N)}=1$. Therefore we have proved the existence of the maximizer for $\nu_{N,s,\alpha,rad}(\mathbb{R}^N)$.
\end{proof}

In the case   $s=0$, applying Lemma 5.1 on Proposition 5.4 we could get the following Corollary.\\
\\
\textbf{Corollary 5.6} \emph{Assume (\ref{equ1.1}) with $s=0$, then
\begin{displaymath}
\nu_{N,0,t,\alpha}(\mathbb{R}^N)=\sup_{\substack{u\in H^{1,N}(\mathbb{R}^N)\\\|\nabla u\|_{L^N(\mathbb{R}^N)}=1}}
\frac{\int_{\mathbb{R}^N}e^{\alpha|u|^{N'}}|u|^N\frac{dx}{|x|^t}}{\|u\|^{N-t}_{L^N(\mathbb{R}^N)}}
\end{displaymath}
is attained.}\\

Then in quite the same way as we prove part (i) of Theorem 1.4, we can prove $\nu_{N,s,t,\alpha}(\mathbb{R}^N)$ is attained by applying Corollary 5.6.\\
\\
{\bf Acknowledgement:} The results of this work has been presented by the first author at an invited talk at the AMS special session 
on Geometric Inequalities and Nonlinear Partial Differential Equations in Las Vegas in April, 2015.

\end{document}